# Towards a better list of citation superstars: compiling a multidisciplinary list of highly cited researchers

**Igor Podlubny**[(1)], **Katarina Kassayova**[(2)]

[(1)]*Department of Applied Informatics and Process Control, BERG Faculty, Technical University of Kosice, Kosice, Slovak Republic*
[(2)]*Department of Physiology, P. J. Safarik University, Kosice, Slovak Republic*

**Abstract:** A new approach to producing multidisciplinary lists of highly cited researchers is described and used for compiling the first multidisciplinary list of highly cited researchers. This approach is essentially related to the recently discovered law of the constant ratios (Podlubny, 2004) and gives a better-balanced representation of different scientific fields.

## Introduction

Citation analysis has become a tool used world-wide for evaluating scientific performance and scientific impact of various subjects operating in research – countries, universities, research institutes, and also individual researchers. The total number of citations is usually considered as an indicator of the scientific impact of the unit under consideration and of the importance of its contribution to the corresponding field of science.

Counting citations has become now, in the era of powerful computers and information technology tools, relatively easy. This provides for administrators of science, management of universities, grant agencies, and others, a seemingly simple tool for justifying their decisions on career promotions, funding and providing other kinds of support of research.

Citation counts are attractive from the viewpoint of presentation of top scientists to the general public. While it is usually difficult to explain to the general public the real scientific contribution of a given researcher and its importance for the advancement of a particular scientific field, the presentation of citation counts is much easier and appears to the public as an objective evaluation.[1] This approach results in so-called lists of highly cited researchers.[2,3] However, those lists are compiled separately for different scientific fields. This is understandable, since comparing, say, ten mathematicians (or ten physicists, or ten chemists, or ten medical researchers, or ten engineers, etc.) using their total citation counts immediately gives a key for sorting them in descending order with respect to their citation counts.

A more difficult problem arises when one has to compare scientists working in different fields, for example, mathematicians, biochemists, physicists, and engineers. The absence of a suitable solution is emphasized by the fact that even the most prolific author of citation analysis, Dr. E. Garfield, used only absolute figures for compiling lists of scientists with the highest impact – see, for example, the list in Ref. 4, where we cannot see any mathematician, engineer, or a specialist in social sciences. The same approach (total numbers of citations) is used also in Ref. 5 entitled "Twenty Years of Citation Superstars", where one can see only specialists in clinical medicine and biomedicine, which does not seem to be adequate with respect to the role of other scientific fields (like mathematics, physics, engineering, etc.) and their contribution to the advancement of science.

The problem of comparing scientific performance in different fields of science is not new. However, previous studies were focused mainly on measuring and comparing the impact of journals or research groups in different fields of subfields of science. Seglen (Ref. 6, Table 1) showed that journal impact factors depend on the research field. Several authors have suggested alternative journal impact measures that account for differences in referencing practices among scientific disciplines. These vary from field-specific impact factors[7] to normalized or relative measures. The latter in principle could enable cross-comparisons of journals among disciplines. For an overview of these efforts see, for example, Refs. 8 and 9.

Similarly, the idea of using normalization for developing more adequate scientometric indicators is also well known. For example, van Raan provided an extensive and detailed studies on field-normalized indicators[10, 11,12]. On the example of chemistry and medicine he demonstrated[11], that field-normalized indicators, such as CPP/FCSm, are different in different fields. However, his study was focused on units of different levels of aggregation, such as research groups and institutes.

An important approach to comparison of different units is based on scaling (or self-similarity) laws observed in scientometrics. Katz[13,14,15] studied scaling relationships between citations counts and the number of publications for research fields, institutes and countries, and suggested scale-independent indicators[15].

Recently, Hirsch[16] suggested a new indicator of citation impact, which is called h-index. He also mentioned that the values of h-index for recognized scientists in different fields are different.

All approaches mentioned above do not deal with comparisons of citation impact of individual scientists working in different scientific fields. In this article we suggest a new approach to compiling multidisciplinary lists of highly cited researchers. This approach is essentially related to the recently discovered law of the constant ratios[17] and gives a better-balanced representation of different scientific fields.

## The data

The main data source for citation counts of highly cited researchers was the ISI Essential Science Indicators (ESI) produced by Thomson ISI.[18] Among many other interesting features, the ESI provides the lists of most frequently cited scientists in twenty-two scientific fields defined by specialists at Thomson ISI for this data source. The citation counts in ESI correspond to ten years plus 2 months period. The data are updated every two months. We used the data from the July 1, 2005 release.

The other source of data was the In-Cites website[3], which, among other interesting information, provides the lists of top ten most frequently cited researchers in the same twenty-two scientific fields. The In-Cites web site is also regularly updated and at each moment covers the moving time window of the same length, namely ten years plus 2 months.

We also used the data from the recent publication of the National Science Foundation, in which the distribution of scientific citations of the U.S. scientific and engineering articles across wide fields of science in 1992, 1994, 1996, 1997, 1999, and 2001 was published (see Ref. 19, Chapter 5, Table 5-27 on page 5-50). The sources for the data appearing in that table were the Science Citation Index (SCI) and the Social Sciences Citation Index (SSCI). Differently from the two aforementioned sources, the NSF report uses nine broad scientific fields. To relate both classifications, we used the definitions of those nine broad fields from the Appendix Table 5-34 on page A5-63 of the NSF report, on one side, and the classification published at the In-Cites website[3], on the other.

## The law of the constant ratio for broad scientific fields

We would like to recall the recently discovered law of the constant ratio[17]: *the ratio of the total number of citations in any two fields of science remains close to constant.* This allows normalization of all fields with respect to mathematics, where the total number of citations is always the smallest, although also growing with time. The law of the constant ratio is in agreement with Katz's studies on scaling and self-similarity in science and scientometric indicators[13,14,15]

In terms of normalized figures, the main result of Ref. 17 can be expressed by Table 1, which says that, for example, one citation in mathematics roughly corresponds to 15 citations in chemistry, 19 citations in physics, and 78 citations in clinical medicine.

**Table 1. Average ratio of the total citation number to the total number of citations in mathematics (see Ref. 17)**

| Field | Average ratio of the total citation number to the total number of citations in mathematics |
|---|---|
| Clinical medicine | 78 |
| Biomedical research | 78 |
| Biology | 8 |
| Chemistry | 15 |
| Physics | 19 |
| Earth/space sciences | 9 |
| Engineering/technology | 5 |
| Mathematics | 1 |
| Social/behavioral sciences | 13 |

We also recall the remark made in Ref. 17, that in the case of top-cited researchers these ratios “will probably need some correction, since the ratios of the peaks in different fields of science do not necessarily copy the ratios shown in Table 1”. Indeed, the total number of citations in a particular scientific field can be interpreted as an integral, while the number of citations of a most frequently cited scientist in that field represents the largest value of the function under integration. Therefore, the question is: does there exist a similar law of constant ratios in the case of highly cited researchers?

## The law of the constant ratio for highly cited researchers

Using the data from In-Cites, we arrived at the conclusion that there *does* exist a similar law, and that there is even certain relationship between the total number of citations in a particular field and the number of citations of the top-cited researcher in that particular field.

Unfortunately, the data available in In-Cites on highly cited researchers do not cover the ten years period – there are only figures for 2003–2005. Those data are updated each month and each of them covers the moving time window of ten years plus two months, so for each of the 22 fields used by In-Cites and ESI we used a set of 24 monthly lists of highly cited researchers.

As the first step, we mapped the ESI/In-Cites classification of 22 scientific fields onto the NSF classification containing 9 fields. For this mapping we used the detailed definitions of the nine broad fields from the Appendix Table 5-34 on page A5-63 of the NSF report[19] and the definitions of the 22 fields published at the In-Cites website[3]. The first two columns of Table 2 represent this mapping. The third column contains the corresponding numbers from Table 1.

Then we performed addition of the ten-element vectors corresponding to those ESI fields that are mapped to the same NSF field (like ESI fields Agriculture and Plant and animal sciences, which are mapped to the NSF field Biology), and subsequently computed the average values of the elements of the resulting vectors. The numbers obtained in this way were considered as the numbers of citations of most frequently cited scientists in the NSF classification fields. These numbers were – similarly to Ref. 17 – normalized with respect to mathematics and are summarized in Table 3. They also appear in the fourth column of Table 2.

The relationship between the average ratio of the total citation number to the total number of citations in mathematics (T), on one side, and the average ratio of the citation counts of ten highly cited scientist in the field to the average number of citations of ten highly cited mathematicians (H), on the other, can be well described by the power law function $H = T^{\alpha}$, where $\alpha=0.82$. If we exclude clinical medicine and biomedicine, a suitable approximation is given even by a much simpler expression: $H = 2T/3$.

**Table 2. Mapping of the ESI fields to NSF fields.**

| ESI field | Broad field according to NSF 2004 report | Average ratio of the total citation number to the total number of citations in mathematics | Average ratio of citations of ten highly cited scientist in the field to the average number of citations of ten highly cited mathematicians |
|---|---|---|---|
| Agriculture | Biology | 8 | 5 |
| Biology and biochemistry | Biomedical research | 78 | 37 |
| Chemistry | Chemistry | 15 | 10 |
| Clinical medicine | Clinical medicine | 78 | 37 |
| Computer science | Engineering and technology | 5 | 3 |
| Economics and business | Social/behavioral sciences | 13 | 9 |
| Engineering | Engineering and technology | 5 | 3 |
| Environment and ecology | Earth and space sciences | 9 | 6 |
| Geo sciences | Earth and space sciences | 9 | 6 |
| Immunology | Clinical medicine | 78 | 37 |
| Material science | Engineering and technology | 5 | 3 |
| Mathematics | Mathematics | 1 | 1 |
| Microbiology | Biomedical research | 78 | 37 |
| Molecular biology and genetics | Biomedical research | 78 | 37 |
| Multidisciplinary | Engineering and technology | 5 | 3 |
| Neuroscience and behavior science | Clinical medicine | 78 | 37 |
| Pharmacology and toxicology | Clinical medicine | 78 | 37 |
| Physics | Physics | 19 | 12 |
| Plant and animal sciences | Biology | 8 | 5 |
| Psychiatry and psychology | Clinical medicine | 78 | 37 |
| Social sciences | Social/behavioral sciences | 13 | 9 |
| Space sciences | Earth and space sciences | 9 | 6 |

**Table 3. The law of the constant ratios for scientific fields and for highly cited researchers in those fields.**

| Average ratio of the total citation number to the total number of citations in mathematics (T) | Average ratio of the citation counts of ten highly cited scientist in the field to the average number of citations of ten highly cited mathematicians (H) |
|---|---|
| 1 | 1 |
| 5 | 3 |
| 8 | 5 |
| 9 | 6 |
| 13 | 9 |
| 15 | 10 |
| 19 | 12 |
| 78 | 37 |

## A multidisciplinary list of highly cited researchers

Using Table 2, we compiled a multidisciplinary list of highly cited researchers. We took the first one hundred of highly cited researchers in each of the twenty-two ESI fields and normalized their citation counts using the fourth column of Table 2. These lists were combined into one, which was then sorted with respect to the number of normalized citations.

The resulting multidisciplinary list of the top 200 highly cited researchers is given in the Appendix to this article. One can see that among the top ten researchers there are five leaders in their particular fields (namely, in Material science, Space science, Chemistry, Mathematics, and Physics), among the top fifty there are leaders in seven fields (add Plant and animal science, and Engineering), and among the top one hundred there are leaders in ten different fields (add Geosciences, Environment and ecology, and Clinical medicine). It could not be so if one used only absolute numbers of citations, like in Ref. 5. In addition, we would like to mention that all twenty-two leaders in their fields are present in the first 1500 lines of the compiled list.

## Discussion

We are far from considering the compiled list as a perfect one.

First, in many cases the numbers of normalized citations of two or more scientists from different fields are very close, so it is difficult to say if their ordering is exact. Instead of considering the exact ordering, it would probably be better to speak about clusters (or groups) of researchers with approximately equivalent citation impact.

Second, one can observe in this list the known problem of possible aggregation of several authors with the same last name and initials (for example, Kobayashi, Nakamura,

Wang, Zhang, etc.). Such cases of possible aggregation are probably indicated by extremely large numbers of articles. For example, around 1000 published articles during 10 years would mean 100 articles during one year, or one published article in 3 or 4 days, which does not seem too realistic. However, according to the correspondence with the Thomson ISI technical support staff [20] and to Ref. 5, splitting such aggregated numbers in parts corresponding to separate persons is currently impossible. To solve this problem, we support the idea of introducing Uniform Author Identifier (UAI)[21], which then would be complementary to the existing and widely used DOI (Document Object Identifier)[22]. If UAI is introduced, it will solve not just the problem of persons with the same names, but also the change of name (especially women's names before and after marriage), different transliterations of non-Latin alphabet names (Russian, Chinese, etc.), misprints in names, and some other situations. Tracking publications and citations using UAI would be much easier and more reliable.

Third, this list is based only on the articles covered by ISI. It does not take into account citations of books, former Soviet journals, many Asian journals, and so forth. As a matter of fact, during the considered period of 10 years Leonhard Euler (EULER L) had more than 530 citations to his original works published in the XVIII-th century. This fact can be easily checked using the Science Citation Index (SCI). This would be sufficient for Euler to be among top 30 mathematicians or among top 400 in our multidisciplinary list, if we do not take into account similar situation of many other past and current authors. However, because Euler's old articles are not in the current ISI database, their citations are not counted.

Fourth, the compiled list does not take into account cross-field citations; all citation numbers are taken as pure citations within the fields. However, cross-field citations are quite common, cross-field citing between almost any pair of fields is asymmetric, and we are not aware of any suitable method for taking this important aspect into account in citation analysis. This topic requires further investigation. Along this way, a new classification scheme of scientific fields and subfields suggested by W. Glänzel and A. Schubert[23] can be used.

Fifth, citation counts, which we used for this study, do not depend on the number of authors of cited papers. However, there is a big difference between an article authored by one author and the article where the list of co-authors contains hundreds of names. For example, the list of authors of Ref. 24 contains 550 names.

Sixth, the ESI and In-Cites are in fact 'black box' products of Thomson Scientific/ISI, so we could not work with the original raw citation data. In our work we assumed that they are based on the same set of raw data as Web of Science. Additional verification of this assumption could be useful.

Our last remark is not related to the numbers that we used for compiling the presented list. It deals with the auxiliary information presented in the last column, namely average number of citations per paper (CPP). In our opinion, it would be better to introduce a new indicator called the average number of citations per "meaningful" paper (CPMP), and set the citation thresholds for "meaningful" papers in different fields. In our list we see a notable number of top-cited researchers who authored one, two, or three articles. It is clear that those who wrote tens or hundreds papers do not have equal response to all of them – some of their articles are highly cited, some other are cited seldom. Therefore, the CPMP could be helpful in balancing these two extreme approaches to producing scientific publications. The number of "meaningful" papers need not to be a fixed number – the h-index, suggested recently by Hirsch[16], is, in fact, a field-independent tool for determining the number of meaningful papers for each author.

## Conclusion

In spite of the above remarks, the presented multidisciplinary list of top cited researchers provides a better picture than Refs. 4 and 5. The approach based on the law of the constant ratio looks like a suitable tool for normalizing citation counts in different fields not only in the case of total numbers, but also in the case of highly cited researchers. Its possible enhancements should preferably solve the problems of equal names of different persons, cross-field citation impact, different numbers of authors of cited articles, and citations of sources that are not in the current ISI database.

## Acknowledgments

The authors are grateful to Milena Matasovska-Tetrevova, the librarian of the University Library of the Technical University of Kosice, for her technical assistance and help, to Dr. Ladislav Pivka for his valuable comments during the preparation of this article, and to the anonymous referees for their suggestions.

# Appendix: Multidisciplinary list of 200 most frequently cited researchers

| | Name | Normalized citations | Field | Rating in the field | Papers | Citations (total) | Citations per paper |
|---|---|---|---|---|---|---|---|
| 1 | INOUE, A | 2495 | MATERIALS SCIENCE | 1 | 655 | 8315 | 12.69 |
| 2 | FILIPPENKO, AV | 1799 | SPACE SCIENCE | 1 | 211 | 10795 | 51.16 |
| 3 | FRENK, CS | 1506 | SPACE SCIENCE | 2 | 123 | 9036 | 73.46 |
| 4 | SCHNEIDER, DP | 1352 | SPACE SCIENCE | 3 | 218 | 8112 | 37.21 |
| 5 | WHITESIDES, GM | 1340 | CHEMISTRY | 1 | 267 | 13399 | 50.18 |
| 6 | ELLIS, RS | 1328 | SPACE SCIENCE | 4 | 127 | 7961 | 62.69 |
| 7 | RAFTERY, AE | 1322 | MATHEMATICS | 1 | 31 | 1322 | 42.65 |
| 8 | NAKAMURA, K | 1292 | PHYSICS | 1 | 549 | 16359 | 29.8 |
| 9 | GRUBBS, RH | 1286 | CHEMISTRY | 2 | 169 | 12853 | 76.05 |
| 10 | LANGDON, TG | 1260 | MATERIALS SCIENCE | 2 | 199 | 4201 | 21.11 |
| 11 | FABIAN, AC | 1251 | SPACE SCIENCE | 5 | 286 | 7505 | 26.24 |
| 12 | WHITE, SDM | 1245 | SPACE SCIENCE | 6 | 109 | 7466 | 68.5 |
| 13 | TOKURA, Y | 1238 | PHYSICS | 2 | 509 | 15666 | 30.78 |
| 14 | STUCKY, GD | 1221 | MATERIALS SCIENCE | 3 | 86 | 4072 | 47.35 |
| 15 | KOBAYASHI, T | 1181 | PHYSICS | 3 | 1000 | 14945 | 14.95 |
| 16 | NAKAMURA, T | 1149 | MATERIALS SCIENCE | 4 | 457 | 3831 | 8.38 |
| 17 | DONOHO, DL | 1137 | MATHEMATICS | 2 | 19 | 1137 | 59.84 |
| 18 | STODDART, JF | 1133 | CHEMISTRY | 3 | 292 | 11331 | 38.8 |
| 19 | KOBAYASHI, S | 1109 | CHEMISTRY | 4 | 679 | 11080 | 16.32 |
| 20 | YORK, DG | 1101 | SPACE SCIENCE | 7 | 174 | 6604 | 37.95 |
| 21 | JOHNSTONE, IM | 1074 | MATHEMATICS | 3 | 14 | 1074 | 76.71 |
| 22 | KASS, RE | 1072 | MATHEMATICS | 4 | 12 | 1072 | 89.33 |
| 23 | MATYJASZEWSKI, K | 1053 | CHEMISTRY | 5 | 278 | 10536 | 37.9 |
| 24 | GUNN, JE | 1047 | SPACE SCIENCE | 8 | 114 | 6284 | 55.12 |
| 25 | XIA, YN | 1046 | MATERIALS SCIENCE | 5 | 70 | 3484 | 49.77 |
| 26 | HECKMAN, TM | 1034 | SPACE SCIENCE | 9 | 151 | 6202 | 41.07 |
| 27 | JONES, JDG | 1019 | PLANT & ANIMAL SCIENCE | 1 | 76 | 5432 | 71.47 |
| 28 | FUKUGITA, M | 1014 | SPACE SCIENCE | 10 | 113 | 6081 | 53.81 |
| 29 | STRAUSS, MA | 1005 | SPACE SCIENCE | 11 | 133 | 6031 | 45.35 |
| 30 | HONSCHEID, K | 1002 | PHYSICS | 4 | 306 | 12690 | 41.47 |
| 31 | HORITA, Z | 992 | MATERIALS SCIENCE | 6 | 127 | 3307 | 26.04 |
| 32 | WANG, J | 990 | CHEMISTRY | 6 | 1075 | 9897 | 9.21 |
| 33 | VANMONTAGU, M | 987 | PLANT & ANIMAL SCIENCE | 2 | 182 | 5263 | 28.92 |
| 34 | COUCH, WJ | 983 | SPACE SCIENCE | 12 | 91 | 5892 | 64.75 |
| 35 | WILLIAMS, DJ | 975 | CHEMISTRY | 7 | 571 | 9744 | 17.06 |
| 36 | IVEZIC, Z | 968 | SPACE SCIENCE | 13 | 132 | 5803 | 43.96 |
| 37 | WANG, J | 941 | ENGINEERING | 1 | 606 | 3134 | 5.17 |
| 38 | BRINKMANN, J | 930 | SPACE SCIENCE | 14 | 186 | 5577 | 29.98 |
| 39 | RHEINGOLD, AL | 926 | CHEMISTRY | 8 | 732 | 9261 | 12.65 |
| 40 | TAKAHASHI, T | 918 | PHYSICS | 5 | 807 | 11630 | 14.41 |
| 41 | VANPARADIJS, J | 915 | SPACE SCIENCE | 15 | 171 | 5491 | 32.11 |
| 42 | HALL, P | 905 | MATHEMATICS | 5 | 162 | 905 | 5.59 |
| 43 | KULKARNI, SR | 891 | SPACE SCIENCE | 16 | 169 | 5348 | 31.64 |
| 44 | RUBIN, DB | 889 | MATHEMATICS | 6 | 27 | 889 | 32.93 |
| 45 | WHITESIDES, GM | 885 | MATERIALS SCIENCE | 7 | 73 | 2951 | 40.42 |
| 46 | GRAB, C | 881 | PHYSICS | 6 | 124 | 11156 | 89.97 |
| 47 | STEIDEL, CC | 879 | SPACE SCIENCE | 17 | 68 | 5270 | 77.5 |
| 48 | COLE, S | 876 | SPACE SCIENCE | 18 | 80 | 5258 | 65.72 |
| 49 | TANAKA, K | 867 | CHEMISTRY | 9 | 1055 | 8676 | 8.22 |
| 50 | NEMOTO, M | 866 | MATERIALS SCIENCE | 8 | 82 | 2883 | 35.16 |
| 51 | BUCHWALD, SL | 863 | CHEMISTRY | 10 | 159 | 8625 | 54.25 |
| 52 | SPEK, AL | 852 | CHEMISTRY | 11 | 587 | 8522 | 14.52 |
| 53 | BENJAMINI, Y | 849 | MATHEMATICS | 7 | 9 | 849 | 94.33 |
| 54 | WEINBERG, DH | 848 | SPACE SCIENCE | 19 | 90 | 5085 | 56.5 |
| 55 | WOLF, G | 846 | PHYSICS | 7 | 600 | 10719 | 17.86 |
| 56 | WANG, J | 836 | PHYSICS | 8 | 1025 | 10585 | 10.33 |
| 57 | WATANABE, Y | 836 | PHYSICS | 9 | 542 | 10585 | 19.53 |
| 58 | WRIGHT, EL | 833 | SPACE SCIENCE | 20 | 51 | 4993 | 97.9 |
| 59 | KLESSIG, DF | 830 | PLANT & ANIMAL SCIENCE | 3 | 73 | 4424 | 60.6 |
| 60 | LUPTON, RH | 830 | SPACE SCIENCE | 21 | 91 | 4975 | 54.67 |
| 61 | KIRSHNER, RP | 828 | SPACE SCIENCE | 22 | 78 | 4972 | 63.74 |
| 62 | CRUTZEN, PJ | 827 | GEOSCIENCES | 1 | 134 | 4958 | 37 |
| 63 | MCMAHON, RG | 824 | SPACE SCIENCE | 23 | 121 | 4942 | 40.84 |
| 64 | GREEN, PJ | 823 | MATHEMATICS | 8 | 16 | 823 | 51.44 |
| 65 | DICKINSON, M | 822 | SPACE SCIENCE | 25 | 76 | 4931 | 64.88 |
| 66 | PETERSON, BA | 822 | SPACE SCIENCE | 24 | 124 | 4932 | 39.77 |
| 67 | VALIEV, RZ | 822 | MATERIALS SCIENCE | 9 | 172 | 2742 | 15.94 |
| 68 | KOKUBO, T | 819 | MATERIALS SCIENCE | 10 | 215 | 2731 | 12.7 |
| 69 | KNAPP, GR | 815 | SPACE SCIENCE | 26 | 119 | 4884 | 41.04 |
| 70 | TRIPPE, TG | 815 | PHYSICS | 10 | 134 | 10313 | 76.96 |
| 71 | SCHLEGEL, DJ | 813 | SPACE SCIENCE | 27 | 55 | 4882 | 88.76 |
| 72 | OSTRIKER, JP | 807 | SPACE SCIENCE | 28 | 94 | 4845 | 51.54 |
| 73 | LEBEDEV, A | 803 | PHYSICS | 11 | 471 | 10168 | 21.59 |
| 74 | YAMADA, S | 803 | PHYSICS | 12 | 430 | 10163 | 23.63 |
| 75 | SMAIL, I | 800 | SPACE SCIENCE | 29 | 119 | 4797 | 40.31 |
| 76 | SUMPTER, JP | 800 | ENVIRONMENT/ECOLOGY | 1 | 38 | 4260 | 112.11 |
| 77 | SUZUKI, T | 795 | PHYSICS | 13 | 1428 | 10070 | 7.05 |
| 78 | CASO, C | 794 | PHYSICS | 14 | 172 | 10054 | 58.45 |
| 79 | CHEONG, SW | 794 | PHYSICS | 15 | 216 | 10045 | 46.5 |
| 80 | EVANS, AG | 794 | MATERIALS SCIENCE | 11 | 129 | 2644 | 20.5 |
| 81 | HERNQUIST, L | 794 | SPACE SCIENCE | 30 | 140 | 4758 | 33.99 |
| 82 | ZHANG, J | 794 | PHYSICS | 16 | 1223 | 10045 | 8.21 |
| 83 | MURAYAMA, H | 789 | PHYSICS | 17 | 145 | 9993 | 68.92 |
| 84 | STAMPFER, MJ | 788 | CLINICAL MEDICINE | 1 | 376 | 30739 | 81.75 |
| 85 | COWIE, LL | 785 | SPACE SCIENCE | 31 | 81 | 4703 | 58.06 |
| 86 | SUZUKI, Y | 785 | PHYSICS | 18 | 537 | 9935 | 18.5 |
| 87 | DIXON, RA | 782 | PLANT & ANIMAL SCIENCE | 4 | 90 | 4170 | 46.33 |
| 88 | CLORE, GM | 777 | CHEMISTRY | 12 | 39 | 7765 | 199.1 |
| 89 | ISHII, T | 777 | PHYSICS | 19 | 316 | 9834 | 31.12 |
| 90 | BAHCALL, NA | 776 | SPACE SCIENCE | 33 | 93 | 4654 | 50.04 |
| 91 | DRESSLER, A | 776 | SPACE SCIENCE | 32 | 61 | 4656 | 76.33 |
| 92 | KOUVELIOTOU, C | 776 | SPACE SCIENCE | 34 | 193 | 4652 | 24.1 |
| 93 | NAVARRO, JF | 774 | SPACE SCIENCE | 35 | 45 | 4640 | 103.11 |
| 94 | AGUILARBENITEZ, M | 771 | PHYSICS | 20 | 206 | 9772 | 47.44 |
| 95 | GURTU, A | 771 | PHYSICS | 21 | 209 | 9757 | 46.68 |
| 96 | NICOLAOU, KC | 771 | CHEMISTRY | 13 | 237 | 7715 | 32.55 |
| 97 | SUNTZEFF, NB | 768 | SPACE SCIENCE | 36 | 80 | 4609 | 57.61 |
| 98 | YAMAMOTO, Y | 765 | CHEMISTRY | 14 | 781 | 7648 | 9.79 |
| 99 | FAN, JQ | 761 | MATHEMATICS | 9 | 46 | 761 | 16.54 |
| 100 | SHINKAI, S | 761 | CHEMISTRY | 15 | 410 | 7601 | 18.54 |
| 101 | HOCHBERG, Y | 760 | MATHEMATICS | 10 | 8 | 760 | 95 |
| 102 | HIGASHI, K | 756 | MATERIALS SCIENCE | 12 | 232 | 2522 | 10.87 |
| 103 | ALIVISATOS, AP | 755 | CHEMISTRY | 16 | 60 | 7548 | 125.8 |
| 104 | EIDELMAN, S | 753 | PHYSICS | 22 | 123 | 9547 | 77.62 |
| 105 | GLAZEBROOK, K | 752 | SPACE SCIENCE | 37 | 88 | 4509 | 51.24 |
| 106 | YAMAGUCHI, K | 752 | CHEMISTRY | 17 | 701 | 7511 | 10.71 |
| 107 | DOI, M | 750 | SPACE SCIENCE | 38 | 92 | 4502 | 48.93 |
| 108 | MARUYAMA, T | 743 | PHYSICS | 23 | 291 | 9410 | 32.34 |
| 109 | REINHOUDT, DN | 741 | CHEMISTRY | 18 | 322 | 7404 | 22.99 |
| 110 | SPERGEL, DN | 740 | SPACE SCIENCE | 39 | 60 | 4440 | 74 |
| 111 | INZE, D | 738 | PLANT & ANIMAL SCIENCE | 5 | 129 | 3938 | 30.53 |
| 112 | BENNETT, CL | 737 | SPACE SCIENCE | 40 | 41 | 4418 | 107.76 |
| 113 | ZHANG, L | 737 | PHYSICS | 24 | 908 | 9326 | 10.27 |
| 114 | MOULD, JR | 735 | SPACE SCIENCE | 41 | 99 | 4409 | 44.54 |
| 115 | BURCHAT, PR | 734 | PHYSICS | 25 | 152 | 9284 | 61.08 |
| 116 | PETTINI, M | 734 | SPACE SCIENCE | 42 | 83 | 4404 | 53.06 |
| 117 | LI, J | 732 | CHEMISTRY | 19 | 954 | 7317 | 7.67 |
| 118 | SZALAY, AS | 732 | SPACE SCIENCE | 43 | 129 | 4392 | 34.05 |
| 119 | HASEGAWA, T | 731 | PHYSICS | 26 | 340 | 9254 | 27.22 |
| 120 | KIM, SB | 731 | PHYSICS | 27 | 268 | 9252 | 34.52 |
| 121 | TILMAN, D | 729 | ENVIRONMENT/ECOLOGY | 2 | 57 | 3887 | 68.19 |
| 122 | AGARWAL, RP | 728 | MATHEMATICS | 11 | 187 | 728 | 3.89 |
| 123 | MUSHOTZKY, RF | 728 | SPACE SCIENCE | 44 | 112 | 4361 | 38.94 |
| 124 | STUCKY, GD | 728 | CHEMISTRY | 20 | 149 | 7278 | 48.85 |
| 125 | SUZUKI, T | 726 | CHEMISTRY | 21 | 1012 | 7254 | 7.17 |
| 126 | FRONTERA, F | 725 | SPACE SCIENCE | 45 | 143 | 4351 | 30.43 |
| 127 | EFSTATHIOU, G | 723 | SPACE SCIENCE | 46 | 102 | 4337 | 42.52 |
| 128 | HUCHRA, JP | 722 | SPACE SCIENCE | 47 | 101 | 4325 | 42.82 |
| 129 | LAMB, DQ | 722 | SPACE SCIENCE | 48 | 113 | 4325 | 38.27 |
| 130 | SHINOZAKI, K | 722 | PLANT & ANIMAL SCIENCE | 6 | 119 | 3845 | 32.31 |
| 131 | FABER, SM | 720 | SPACE SCIENCE | 49 | 70 | 4320 | 61.71 |
| 132 | FRAIL, DA | 720 | SPACE SCIENCE | 50 | 126 | 4319 | 34.28 |
| 133 | FRECHET, JMJ | 720 | CHEMISTRY | 22 | 222 | 7204 | 32.45 |
| 134 | FENG, JL | 719 | PHYSICS | 28 | 68 | 9102 | 133.85 |
| 135 | MIRKIN, CA | 719 | CHEMISTRY | 23 | 130 | 7179 | 55.22 |
| 136 | MUNN, JA | 719 | SPACE SCIENCE | 51 | 75 | 4314 | 57.52 |
| 137 | CHIB, S | 717 | MATHEMATICS | 12 | 10 | 717 | 71.7 |
| 138 | CORMA, A | 716 | CHEMISTRY | 24 | 324 | 7159 | 22.1 |
| 139 | TROST, BM | 716 | CHEMISTRY | 25 | 231 | 7158 | 30.99 |
| 140 | CSABAI, I | 714 | SPACE SCIENCE | 52 | 90 | 4287 | 47.63 |
| 141 | GENZEL, R | 714 | SPACE SCIENCE | 53 | 132 | 4281 | 32.43 |
| 142 | FAN, XH | 713 | SPACE SCIENCE | 54 | 81 | 4279 | 52.83 |
| 143 | NAKAMURA, S | 711 | PHYSICS | 29 | 391 | 9006 | 23.03 |
| 144 | BAULCOMBE, DC | 707 | PLANT & ANIMAL SCIENCE | 7 | 45 | 3765 | 83.67 |
| 145 | HERRMANN, WA | 707 | CHEMISTRY | 26 | 195 | 7058 | 36.19 |
| 146 | OZIN, GA | 701 | MATERIALS SCIENCE | 13 | 95 | 2336 | 24.59 |
| 147 | SUTHERLAND, W | 701 | SPACE SCIENCE | 55 | 101 | 4205 | 41.63 |
| 148 | AMSLER, C | 699 | PHYSICS | 30 | 73 | 8845 | 121.16 |
| 149 | HINSHAW, G | 698 | SPACE SCIENCE | 56 | 40 | 4187 | 104.67 |
| 150 | KALNAY, E | 696 | GEOSCIENCES | 2 | 29 | 4172 | 143.86 |
| 151 | TIELENS, AGGM | 696 | SPACE SCIENCE | 57 | 161 | 4177 | 25.94 |
| 152 | RIESS, AG | 693 | SPACE SCIENCE | 58 | 41 | 4157 | 101.39 |
| 153 | TIBSHIRANI, R | 692 | MATHEMATICS | 13 | 28 | 692 | 24.71 |
| 154 | BRANDT, WN | 692 | SPACE SCIENCE | 60 | 154 | 4147 | 26.93 |
| 155 | KOGUT, A | 692 | SPACE SCIENCE | 59 | 48 | 4152 | 86.5 |
| 156 | WANG, Y | 692 | CHEMISTRY | 27 | 1326 | 6915 | 5.21 |
| 157 | ZHANG, Y | 692 | MATERIALS SCIENCE | 14 | 504 | 2304 | 4.57 |
| 158 | TIERNEY, L | 690 | MATHEMATICS | 14 | 8 | 690 | 86.25 |
| 159 | REYNOLDS, RW | 689 | GEOSCIENCES | 3 | 22 | 4134 | 187.91 |
| 160 | GIAVALISCO, M | 687 | SPACE SCIENCE | 61 | 46 | 4118 | 89.52 |
| 161 | KIM, HJ | 687 | PHYSICS | 31 | 743 | 8697 | 11.71 |
| 162 | NOYORI, R | 687 | CHEMISTRY | 28 | 110 | 6875 | 62.5 |
| 163 | ZHANG, T | 684 | MATERIALS SCIENCE | 15 | 175 | 2281 | 13.03 |
| 164 | READ, RJ | 683 | CHEMISTRY | 29 | 17 | 6832 | 401.88 |
| 165 | PHILLIPS, MM | 680 | SPACE SCIENCE | 62 | 54 | 4073 | 75.43 |
| 166 | ITO, Y | 678 | CHEMISTRY | 30 | 753 | 6782 | 9.01 |
| 167 | SIMONSON, T | 678 | CHEMISTRY | 31 | 20 | 6781 | 339.05 |
| 168 | SEEBACH, D | 677 | CHEMISTRY | 32 | 181 | 6772 | 37.41 |
| 169 | VITOUSEK, PM | 677 | ENVIRONMENT/ECOLOGY | 3 | 75 | 3605 | 48.07 |
| 170 | WU, X | 675 | PHYSICS | 32 | 539 | 8547 | 15.86 |
| 171 | BRUNGER, AT | 674 | CHEMISTRY | 33 | 17 | 6738 | 396.35 |
| 172 | ANNIS, J | 672 | SPACE SCIENCE | 63 | 65 | 4031 | 62.02 |
| 173 | KUSZEWSKI, J | 671 | CHEMISTRY | 34 | 12 | 6708 | 559 |
| 174 | LI, J | 671 | PHYSICS | 33 | 971 | 8492 | 8.75 |
| 175 | KANAMITSU, M | 669 | GEOSCIENCES | 4 | 26 | 4010 | 154.23 |
| 176 | FREEMAN, KC | 668 | SPACE SCIENCE | 64 | 186 | 4008 | 21.55 |
| 177 | YANG, PD | 668 | MATERIALS SCIENCE | 16 | 27 | 2224 | 82.37 |
| 178 | KOLLMAN, PA | 666 | CHEMISTRY | 36 | 107 | 6654 | 62.19 |
| 179 | OLIVE, KA | 666 | PHYSICS | 34 | 84 | 8429 | 100.35 |
| 180 | PANNU, NS | 666 | CHEMISTRY | 35 | 9 | 6657 | 739.67 |
| 181 | GROS, P | 665 | CHEMISTRY | 37 | 53 | 6647 | 125.42 |
| 182 | LEETMAA, A | 665 | GEOSCIENCES | 5 | 21 | 3986 | 189.81 |
| 183 | MIKOS, AG | 665 | MATERIALS SCIENCE | 17 | 77 | 2216 | 28.78 |
| 184 | PIER, JR | 665 | SPACE SCIENCE | 65 | 64 | 3988 | 62.31 |
| 185 | SMALLEY, RE | 663 | PHYSICS | 35 | 58 | 8401 | 144.84 |
| 186 | WHITE, AJP | 663 | CHEMISTRY | 38 | 402 | 6635 | 16.5 |
| 187 | ILLINGWORTH, GD | 662 | SPACE SCIENCE | 66 | 96 | 3968 | 41.33 |
| 188 | CRANDALL, KA | 660 | COMPUTER SCIENCE | 1 | 2 | 2201 | 1100.5 |
| 189 | MO, KC | 660 | GEOSCIENCES | 6 | 35 | 3961 | 113.17 |
| 190 | POSADA, D | 660 | COMPUTER SCIENCE | 2 | 3 | 2199 | 733 |
| 191 | ZHANG, Y | 660 | PHYSICS | 36 | 1022 | 8357 | 8.18 |
| 192 | DAVIS, M | 659 | SPACE SCIENCE | 68 | 54 | 3947 | 73.09 |
| 193 | FOYER, CH | 659 | PLANT & ANIMAL SCIENCE | 8 | 99 | 3512 | 35.47 |
| 194 | HO, LC | 659 | SPACE SCIENCE | 67 | 116 | 3955 | 34.09 |
| 195 | YAGHI, OM | 659 | CHEMISTRY | 39 | 57 | 6584 | 115.51 |
| 196 | GRATZEL, M | 654 | CHEMISTRY | 40 | 138 | 6533 | 47.34 |
| 197 | HAGIWARA, K | 654 | PHYSICS | 37 | 68 | 8288 | 121.88 |
| 198 | ROBINS, JM | 653 | MATHEMATICS | 15 | 31 | 653 | 21.06 |
| 199 | KOBAYASHI, K | 653 | PHYSICS | 38 | 446 | 8264 | 18.53 |
| 200 | JACOB, DJ | 651 | GEOSCIENCES | 7 | 127 | 3904 | 30.74 |